\documentclass[12pt]{article}

\usepackage{amsmath,amssymb,amsfonts,amsthm,graphics,verbatim}

\newtheorem{lemma}{Lemma}[subsection]

\newtheorem{thm}{Theorem}[subsection]

\theoremstyle{definition}
\newtheorem*{remark trees}{Tree boundaries}
\newtheorem*{remark tits}{Tits boundaries}
\newtheorem*{remark note}{Setting notation}
\newtheorem*{remark}{Remark}

\DeclareMathOperator{\G}{\bold{G}}

\begin{document}

\newcommand{\ka}{\kappa}
\newcommand{\ga}{\Gamma}
\newcommand{\tga}{\Gamma}
\newcommand{\tg}{H}
\newcommand{\rt}{\rightarrow}
\newcommand{\gmga}{\tga \backslash \tg}
\newcommand{\neu}{N(\ga)}
\newcommand{\ap}{{\mathcal{A}}}
\newcommand{\lb}{\lambda}
\newcommand{\se}{\subseteq}
\newcommand{\e}{\varepsilon}
\newcommand{\nre}{\sqrt[n]{\varepsilon}}
\newcommand{\tp}{\widetilde{\phi}}
\newcommand{\pp}{\partial \phi}
\newcommand{\app}{{\mathcal{A}}^+}
\newcommand{\ppr}{\partial \phi _R}
\newcommand{\pps}{{\partial \phi _R}_*}
\newcommand{\pu}{U _{\partial}}
\newcommand{\gh}{\Gamma}
\newcommand{\ghp}{\Gamma ^{\varphi}}
\newcommand{\hp}{G^+}
\newcommand{\oppr}{\overline{\partial \phi _R}}
\newcommand{\ta}{K \gamma_\mathfrak{p}}
\newcommand{\de}{\delta}
\newcommand{\tgd}{K \gamma_\mathfrak{p}(\de)}
\newcommand{\p}{\phi}
\newcommand{\pn}{\phi ^{-1}}
\newcommand{\g}{\gamma}
\newcommand{\pgd}{\pi(\gr ,\de)}
\newcommand{\ep}{\hfill $\blacksquare$ \bigskip}
\newcommand{\pru}{\bigskip \noindent {\bf Proof:} \;}
\newcommand{\sa}{\measuredangle _}
\newcommand{\xy}{X_\infty}
\newcommand{\gy}{\gamma_\infty}
\newcommand{\py}{\pi_\infty}
\newcommand{\ey}{e_\infty}
\newcommand{\xr}{X_{\mathfrak{p}}}
\newcommand{\gr}{\gamma_{\mathfrak{p}}}
\newcommand{\pr}{\pi_{\mathfrak{p}}}
\newcommand{\er}{e_{\mathfrak{p}}}
\newcommand{\phr}{\phi_{\mathfrak{p}}}
\newcommand{\phy}{\phi_{\infty}}
\newcommand{\ar}{\G (\mathcal{O}_S)}
\newcommand{\Ar}{\G (\mathcal{O}_S)}

\title{Quasi-isometries of rank one $S$-arithmetic lattices}
\author{Kevin Wortman\thanks{Supported in part by an N.S.F. Postdoctoral
Fellowship.}}
\date{December 7, 2007}
\maketitle

\begin{abstract}

We complete the quasi-isometric classification of irreducible lattices in semisimple Lie groups over nondiscrete locally compact fields of characteristic zero by showing that any quasi-isometry of a rank one $S$-arithmetic lattice in a semisimple Lie group over nondiscrete locally compact fields of characteristic zero is a finite distance in the $\sup$-norm from a commensurator.

\end{abstract}

\begin{section}{Introduction}

Throughout we let $K$ be an algebraic number field, $V_K$ the set
of all inequivalent valuations on $K$, and $V_K^\infty \subseteq
V_K$ the subset of archimedean valuations. We will use $S$ to
denote a finite subset of $V_K$ that contains $V_K^\infty$, and we
write the corresponding ring of $S$-integers in $K$ as
$\mathcal{O}_S$. In this paper, $\mathbf{G}$ will always be a connected
non-commutative absolutely simple algebraic $K$-group.

\begin{subsection}{Commensurators} For any valuation
$v\in V_K$, we let $K_v$ be the completion of $K$ with respect to
$v$. For any set of valuations $S' \se V_K$, we define
$$G_{S'} = \prod _{v \in S'} \mathbf{G}(K_v)$$
and we identify $\mathbf{G}(\mathcal{O}_S)$
 as a discrete subgroup of $G_S $ using the diagonal embedding.

We let $\text{Aut}(G_S)$ be the group of topological group
automorphisms of $G_S$. An automorphism $\psi \in \text{Aut}(G_S)$
\emph{commensurates} $\mathbf{G}(\mathcal{O}_S)$ if
$\psi(\mathbf{G}(\mathcal{O}_S)) \cap \mathbf{G}(\mathcal{O}_S)$
is a finite index subgroup of both $\psi
(\mathbf{G}(\mathcal{O}_S) )$ and $\mathbf{G}(\mathcal{O}_S)$.

We define the \emph{commensurator group of}
$\mathbf{G}(\mathcal{O}_S)$ to be the subgroup of
$\text{Aut}(G_S)$ consisting of automorphisms that commensurate
$\mathbf{G}(\mathcal{O}_S)$. This group is denoted as
$\text{Comm}_{\text{Aut}(G_S)}(\mathbf{G}(\mathcal{O}_S))$. Notice that it
differs from the standard definition of the commensurator group
of $\mathbf{G}(\mathcal{O}_S)$ in that we have not restricted
ourselves to inner automorphisms.
\end{subsection}

\begin{subsection}{Quasi-isometries} For constants $L \geq
1$ and $C \geq 0$, an $(L,C)$ \emph{quasi-isometric embedding} of
a metric space $X$ into a metric space $Y$ is a function $\phi : X
\rightarrow Y$ such that for any $x_{1},x_{2} \in X$:
$$\frac{1}{L} d\big(x_{1},x_{2}\big) -C \leq d\big(\phi(x_{1}),\phi(x_{2})\big)
\leq   L d\big(x_{1},x_{2}\big) +C$$

We call $\phi$ an $(L,C)$ \emph{quasi-isometry} if $\phi$ is an
$(L,C)$ quasi-isometric embedding and there is a number $D\geq 0$
such that every point in $Y$ is within distance $D$ of some point
in the image of $X$.

\end{subsection}

\begin{subsection}{Quasi-isometry groups} For a metric space $X$,
we define the relation $\sim$  on the set of functions  $X\rt X$
by $\phi \sim \psi$ if
$$\sup _{x\in X} d\big(\phi(x),\psi(x)\big)<\infty$$ In this paper we will call
two functions \emph{equivalent} if they are related by $\sim$.

For a finitely generated group with a word metric $\Gamma$, we
form the set of all quasi-isometries of $\Gamma$, and denote the
quotient space modulo $\sim$ by $\mathcal{QI}(\Gamma)$. We call
$\mathcal{QI}(\Gamma)$ the \emph{quasi-isometry group} of $\Gamma$
as it has a natural group structure arising from function
composition.

\end{subsection}

\begin{subsection}{Main result}

In this paper we show:

\begin{thm}\label{t:a} Suppose $K$ is an algebraic number field,
$\mathbf{G}$ is a connected
non-commutative absolutely simple algebraic $K$-group, $S$ properly contains $V_K^\infty$, and $\text{\emph{rank}}_K(\mathbf{G})=1$. Then there is an isomorphism
$$ \mathcal{QI}(\mathbf{G}(\mathcal{O}_S) ) \cong \text{\emph{Comm}}_{\text{\emph{Aut}}(G_S)}(\mathbf{G}(\mathcal{O}_S))$$

\end{thm}

Two special cases of Theorem~\ref{t:a} had been previously known:
Taback proved it when $\mathbf{G}(\mathcal{O}_S)$ is commensurable
to $\mathbf{PGL_2}(\mathbb{Z}[1/p])$ where $p$ is a prime number
\cite{Ta} (we'll use Taback's theorem in our proof), and it was
proved when $\text{rank}_{K_v}(\mathbf{G}) \geq 2$ for all $v \in
S$ by the author in \cite{W2}.

Examples of $S$-arithmetic groups for which Theorem~\ref{t:a} had
been previously unknown include when
$\mathbf{G}(\mathcal{O}_S)=\mathbf{PGL_2}(\mathbb{Z}[1/m])$ where
$m$ is composite. Theorem~\ref{t:a} in this case alone was an
object of study; see Taback-Whyte \cite{T-W} for their program of
study. Theorem~\ref{t:m} below presents a short proof of this
case.

\end{subsection}

\begin{subsection}{Quasi-isometry groups of non-cocompact irreducible $S$-arithmetic lattices} Combining Theorem~\ref{t:a} with the results from Schwartz, Farb-Schwartz, Eskin, Farb, Taback, and Wortman (\cite{S}, \cite{F-S}, \cite{S 2}, \cite{Es}, \cite{Fa}, \cite{Ta}, and \cite{W2}) we have

\begin{thm}\label{t:b} Suppose $K$ is an algebraic number field, and
$\mathbf{G}$ is a connected,
non-commutative, absolutely simple, $K$-isotropic, algebraic $K$-group.
If either $K \ncong \mathbb{Q}$,
$S \neq V_K^\infty$, or $\mathbf{G}$ is not
$\mathbb{Q}$-isomorphic to $\mathbf{PGL_2}$, then there is an
isomorphism
$$ \mathcal{QI}(\mathbf{G}(\mathcal{O}_S) ) \cong \text{\emph{Comm}}_{\text{\emph{Aut}}(G_S)}(\mathbf{G}(\mathcal{O}_S))$$

\end{thm}

Note that Theorem~\ref{t:b} identifies the quasi-isometry group of
any non-cocompact irreducible $S$-arithmetic lattice in a
semisimple Lie group over nondiscrete locally compact fields of
characteristic $0$ that is not virtually free. Indeed, any such
lattice is isomorphic ``up to finite groups" to some
$\mathbf{G}(\mathcal{O}_S)$ fitting the hypothesis of
Theorem~\ref{t:b}.

\end{subsection}

\begin{subsection}{Function fields} Our proof of
Theorem~\ref{t:a} also applies when $K$ is a global function field
when we add the hypothesis that there exists $v,w \in S$ such that
$\text{rank}_{K_v}(\mathbf{G}) =1$ and
$\text{rank}_{K_w}(\mathbf{G}) \geq 2$.

For more on what is known about the quasi-isometry groups of
arithmetic groups over function fields -- and for a conjectural
picture of what is unknown -- see \cite{ff}.

\end{subsection}

\begin{subsection}{Outline of the proof and the paper} We begin
Section~\ref{s:1} with a sort of large-scale reduction theory. We
examine a metric neighborhood, $N$, of an orbit of an
$S$-arithmetic group, $\Gamma$, inside the natural product of
Euclidean buildings and symmetric spaces, $X$. In
Section~\ref{s:2.2} we show that the fibers of $N$ under
projections to building factors of $X$ are geometric models for
$S$-arithmetic subgroups of $\Gamma$.

In Section~\ref{s:2.3}, we apply results from \cite{W2} to extend
a quasi-isometry $\phi: \Gamma \rightarrow \Gamma$ to the space
$X$ that necessarily preserves factors.

Our general approach to proving Theorem~\ref{t:a} is to restrict
$\phi$ to factors of $X$ and use the results from
Section~\ref{s:2.2} to decompose $\phi$ into quasi-isometries of
$S$-arithmetic subgroups of $\Gamma$. Once each of these
``sub-quasi-isometries'' is understood, they are pieced together
to show that $\phi$ is a commensurator. An easy example of this
technique is given in Section~\ref{s:2.4} by
$\mathbf{PGL_2}(\mathbb{Z}[1/m])$. We then treat the general case
of Theorem~\ref{t:a} in Section~\ref{s:2.5}

Our proof in Section~\ref{s:2.5} relies on the structure of
horoballs for $S$-arithmetic groups associated to the product of a
symmetric space and a single tree. We prove the results that we
need for these horoballs (that they are connected, pairwise
disjoint, and reflect a kind of symmetry in factors) in
Section~\ref{s:2}, so that Section~\ref{s:2} is somewhat of an
appendix. The proof is organized in this way because, in the
author's opinion, it just makes it easier to digest the material.
There would be no harm though in reading Section~\ref{s:2} after
Section~\ref{s:2.4} and before Section~\ref{s:2.5} for those who
prefer a more linear presentation.

\end{subsection}

\begin{subsection}{Acknowledgements} I was fortunate to have several conversations with Kevin Whyte on the contents of this paper, and am thankful for those. In particular, he brought to my attention  that $\mathbf{SL_2}(\mathbb{Z})$ is the only non-cocompact arithmetic lattice in $\mathbf{SL_2}(\mathbb{R})$ up to commesurability.

I am also grateful to the following mathematicians
who contributed to this paper: Tara Brendle, Kariane Calta, Indira
Chatterji, Alex Eskin, Benson Farb, Dan Margalit, Dave Witte
Morris, and Jennifer Taback.

\end{subsection}

\end{section}

\begin{section}{Proof of Theorem~\ref{t:a}}\label{s:1}

Let $\mathbf{G}$, $K$, and $S$ be as in Theorem~\ref{t:a}, and let $\phi :
\mathbf{G}(\mathcal{O}_S) \rt \mathbf{G}(\mathcal{O}_S)$ be a
quasi-isometry.

\begin{subsection}{Geometric models}
For each valuation $v$ of $K$, we let $X_v$ be the symmetric space
or Euclidean building corresponding to $\mathbf{G}(K_v)$. If $S'$
is a finite set of valuations of $K$, we let
$$X_{S'}=\prod_{v\in S'}X_v$$ Recall that there is a natural
inclusion of topological groups $\text{Aut} (G_{S'})
\hookrightarrow \text{Isom}(X_{S'})$.

Let $\mathcal{O}$ be the ring of integers in $K$, and fix a
connected subspace $\Omega _ {V_K^\infty} \se X_{V_K^\infty}$ that
$\mathbf{G}(\mathcal{O})$ acts cocompactly on. Let $D_\infty \se
X_{V_K^\infty}$ be a fundamental domain for this action.

For each nonarchimedean valuation $w \in S-V_K^\infty$, we denote
the ring of integers in $K_w$ by $\mathcal{O}_w$. The group
$\mathbf{G}(\mathcal{O}_w)$ is bounded in $\mathbf{G}(K_w)$,
 so $\mathbf{G}(\mathcal{O}_w)$ fixes a point $x _w \in
X_w$.
 We choose a bounded set $D_w \se X_w$ containing $x_w$ with
  $\mathbf{G}(\mathcal{O}_S) D_w =X_w$ and such that
 $g x_w \in D_w$ for $g \in \mathbf{G}(\mathcal{O}_S)$ implies that $gx_w=x_w$.

  For any set of
valuations $S'$ satisfying $V_K^\infty \se S' \se S$, we define
the space
  $$\Omega _{S'} = \mathbf{G}(\mathcal{O}_{S'})\Big(D _\infty \times
  \prod_{w \in S'-V_K^\infty }D _w\Big)$$
Note that $\Omega_{S'}$ is a subspace of $ X_{S'}$.

 We endow
$\Omega_{S'}$ with the path metric. Since
$\mathbf{G}(\mathcal{O}_{S'})$ acts cocompactly on $\Omega_{S'}$,
we have the following observation:

\begin{lemma}\label{l:1.1} For $V_K^\infty \se S' \se S$, the space $\Omega
_{S'}$ is quasi-isometric to the group
$\mathbf{G}(\mathcal{O}_{S'})$.
 \end{lemma}

\end{subsection}

\begin{subsection}{Fibers of projections to buildings are
$S$-arithmetic}\label{s:2.2} In the large-scale, the fibers of the
projection of $\Omega _S$ onto building factors of $X_S$ are also
$S$-arithmetic groups (or more precisely, $S'$-arithmetic groups).
This is the statement of Lemma~\ref{l:1.3} below, but we will
start with a proof of a special case.

\begin{lemma}\label{l:1.2} The Hausdorff distance between
 $$\Omega _S
\cap\Big(X_{S'} \times \prod_{w \in S-S'} \{ x_w \}\Big)$$
and
$$\Omega _{S'} \times \prod _{w \in S-S'}\{x_w\}$$
is finite.

\end{lemma}

\begin{proof} There are three main steps in this
proof.

First, if $y \in \Omega _{S'}$, then $y =gd $ for some $g \in
\mathbf{G}(\mathcal{O}_{S'})$ and some $$d \in D_\infty \times
\prod_{w\in S' -V_K^\infty}D_w$$ Since
$\mathbf{G}(\mathcal{O}_{S'})\leq \mathbf{G}(\mathcal{O}_w)$ for
all $w \in S-S'$, it follows from our choice of the points $x_w$
that
\begin{align*} \{y\} \times \prod_{w\in S-S'}\{x_w\} & = g \Big(\{d\}\times
\prod _{w\in S-S'}\{x_w\}\Big) \se \Omega _S
\end{align*}
Therefore,
\begin{align} \Omega _{S'} \times \prod _{w \in S-S'}\{x_w\}
 \se \Omega _S \cap \Big(X_{S'} \times \prod_{w \in S-S'} \{ x_w \}\Big)
\end{align}

Second, we suppose $$\{y\}\times \prod _{w \in S-S'}\{x_w\} \se
\Omega _S$$ for some $y \in X_{S'}$. Then there exists a $g \in
\mathbf{G}(\mathcal{O}_S)$ such that $$gy \in D_\infty \times
\prod _{w \in S'-V_K^\infty} D_w$$ and $g x_w \in D_w$ for all $w
\in S-S'$. Notice that our choice of $D_w$ implies $g x_w =x_w$
for all $w \in S-S'$. Thus, $g$ is contained in the compact group
$$H_w=\{\,h \in \mathbf{G}(K_w) \mid hx_w=x_w\,\}$$ for all $w
\in S-S'$. Consequently, g is contained in the discrete group
$$\mathbf{G}(\mathcal{O}_S) \cap \Big(G_{S'} \times \prod_{w \in S-S'}H_w \Big)$$
We name this discrete group $\Gamma _{S'}$.

Note that we have shown
$$\{y\} \times \prod _{w \in S-S'}\{x_w\} \se \Gamma _{S'} \Big( D_\infty \times \prod_{w \in S-V_K^\infty}
D_w \Big)$$ Therefore, \begin{align}\Omega _S \cap\Big(X_{S'}
\times \prod_{w \in S-S'} \{ x_w \}\Big) \se \Gamma _{S'} \Big(
D_\infty \times \prod_{w \in S-V_K^\infty} D_w \Big) \end{align}

Third, we recall that
$$\mathbf{G}(\mathcal{O}_{S'}) = \mathbf{G}(\mathcal{O}_S) \cap
\Big(G_{S'} \times \prod_{w \in S-S'}\mathbf{G}(\mathcal{O}_w)
\Big)$$ and use the definition of $\Gamma _{S'}$ coupled with the
fact that $\mathbf{G}(\mathcal{O}_w) \leq H_w$ to see that $\Gamma
_{S'}$ contains $\mathbf{G}(\mathcal{O}_{S'})$. Since, $\Gamma
_{S'}$ and $\mathbf{G}(\mathcal{O}_{S'})$ are lattices in $G_{S'}
\times \prod _{w \in S-S'}H_w$, the containment
$\mathbf{G}(\mathcal{O}_{S'}) \leq \Gamma _{S'}$ is of finite
index. Therefore, the Hausdorff distance between
$$\Gamma _{S'} \Big( D_\infty \times \prod_{w \in S' -V_K^\infty} D_w \times \prod_{w \in S -S'} \{x_w\}
 \Big)$$
and
$$\Omega_{S'} \times \prod_{w\in S-S'} \{x_w \} = \mathbf{G}(\mathcal{O}_{S'}) \Big( D_\infty \times \prod_{w \in S' -V_K^\infty} D_w \times \prod_{w \in S -S'} \{x_w\}
 \Big)$$
is finite. Combined with (1) and (2) above, the lemma follows.

\end{proof}

The more general form of Lemma~\ref{l:1.2} that we will use in our proof
of Theorem~\ref{t:a} is the following lemma. We will use the notation of
$x_{S-S'}$ for the point $(x_w)_{w\in S-S'} \in X_{S-S'}$.

\begin{lemma}\label{l:1.3} Suppose $V_K^\infty \se S' \se S$. If $y \in X_{S-S'}$
and $y \in \mathbf{G}(\mathcal{O}_S) x_{ S-S'}$,
 then the Hausdorff distance between $$\Omega _S \cap\Big(X_{S'} \times \{ y \}\Big)$$
and
 $$ \Omega _{S'}
\times \{ y \}$$ is finite.

\end{lemma}

\begin{remark} Our assumption in Lemma~\ref{l:1.3} that $y \in
\mathbf{G}(\mathcal{O}_S) x_{ S-S'}$ is not a serious restriction
over the assumption that $y \in X_{S-S'}$. Indeed,
$\mathbf{G}(\mathcal{O}_S) $ is dense in $G_{S-S'}$, so the orbit
$\mathbf{G}(\mathcal{O}_S) x_{ S-S'}$ is a finite Hausdorff
distance from the space $X_{ S-S'}$.

\end{remark}

\begin{proof} Let $g \in \mathbf{G}(\mathcal{O}_S)$ be
such that $ y = g x_{ S-S'}$.  Then
\begin{align*}
\{\,h \in \mathbf{G}(\mathcal{O}_S) \mid hx_{S-S'} =y \,\} & = g
\{\,h \in \mathbf{G}(\mathcal{O}_S) \mid  hx_{S-S'} =x_{S-S'} \,\}
\\
 & = g \Big( \mathbf{G}(\mathcal{O}_S) \cap \Big(G_{S'} \times \prod_{w \in S-S'}H_w \Big) \Big)
  \\
 & = g \Gamma _{S'}
 \end{align*}
 where $H_w$ and $\Gamma _{S'}$ are as in the proof of the
 previous lemma.

 Now by our choice of the points $x_w \in
X_w$ for $w \in S-V_K^\infty$ at the beginning of this section, we
have
\begin{align*} \Omega _S \cap \Big(X_{S'} \times \{y \} \Big) & =
\mathbf{G}(\mathcal{O}_S)\Big(D_\infty \times \prod_{w\in S -
V_K^\infty} D_w \Big) \cap \Big(X_{S'} \times \{y \} \Big) \\
& = g \Gamma _{S'} \Big(D_\infty \times \prod_{w\in S' -
V_K^\infty} D_w \Big) \times \{y\}
\end{align*}

Notice that the final space from the above chain of equalities is
a finite Hausdorff distance from $$g \mathbf{G}(\mathcal{O} _{S'})
\Big(D_\infty \times \prod_{w\in S' - V_K^\infty} D_w \Big) \times
\{y\}$$ since $\Gamma _{S'} $ is commensurable with
$\mathbf{G}(\mathcal{O}_{S'})$.

Because $g$ commensurates $\mathbf{G}(\mathcal{O}_{S'})$, the
above space is also a finite Hausdorff distance from $\Omega _{S'}
\times \{y\}$.

\end{proof}

\end{subsection}

\begin{subsection}{Extending quasi-isometries of $\Omega _S$ to
$X_S$}\label{s:2.3} Applying Lemma~\ref{l:1.1}, we can regard our
quasi-isometry $\phi : \mathbf{G}(\mathcal{O}_S) \rt
\mathbf{G}(\mathcal{O}_S) $ as a quasi-isometry of $\Omega _S$.
Our goal is to show that $\phi $ is equivalent to an element of
$\text{Comm}_{\text{Aut}(G_S)}(\mathbf{G}(\mathcal{O}_S))$, and we
begin by extending $\phi$ to a quasi-isometry of $X_S$.

\begin{lemma}\label{l:1.4} There is a permutation of $S$, which we
name $\tau$, and there are quasi-isometries $$\phi _v : X_v
\rightarrow X_{\tau (v)}$$ such that the restriction of the
quasi-isometry
$$\phi _S = \prod _{v \in S} \phi_v : X_S \rt X_S$$ to $\Omega _S$
is equivalent to $\phi$. If $X_v$ is a higher rank space for any $v \in S$, then
$\phi _v$ may be taken to be an isometry.

\end{lemma}

\begin{proof} By Proposition 6.9 of \cite{W2}, the quasi-isometry $\phi :
\Omega _S \rightarrow \Omega _S $ extends to a quasi-isometry of
$X$. That is there is some quasi-isometry
$\overline{\phi}:X\rightarrow X$ such that $\phi \sim
\overline{\phi}|_{\Omega _S}$ where $\overline{\phi}|_{\Omega _S}$
is the restriction of $\overline{\phi}$ to $\Omega _S$.

The map $\overline{\phi}$ preserves factors in the boundary of $X$
and an argument of Eskin's -- Proposition 10.1 of \cite{Es} -- can
be directly applied to show that $\overline{\phi}$ is equivalent
to a product of quasi-isometries of the factors of $X$, up to
permutation of factors.

Note that the statement of Proposition 10.1 from \cite{Es} claims
that $X_v$ and $X_{\tau (v)}$ are isometric for $v \in
V_K^\infty$. This is because quasi-isometric symmetric spaces are
isometric up to scale.

\end{proof}

\end{subsection}

\begin{subsection}{Example of proof to come}\label{s:2.4} Before continuing with the
 general proof, we'll pause for a moment to demonstrate the
  utility of Lemmas~\ref{l:1.3} and~\ref{l:1.4} by proving the following
  special case of Theorem~\ref{t:a}.

\begin{thm}\label{t:m} If $m\in \mathbb{N}$ and $m\neq 1$, then
$$\mathcal{QI}\big(\mathbf{PGL_2}(\mathbb{Z}[1/m])\big)\cong
\mathbf{PGL_2}(\mathbb{Q})$$
\end{thm}

\begin{proof}
Let $K=\mathbb{Q}$, $\mathbf{G}=\mathbf{PGL_2}$, and
$S=\{v_\infty\}\cup\{v_p\}_{p|m}$ where $v_\infty$ is the
archimedean valuation and $v_p$ is the $p$-adic valuation. Thus,
$\mathbf{G}(\mathcal{O}_S)=\mathbf{PGL_2}(\mathbb{Z}[1/m])$, the
space $X_{v_\infty}$ is the hyperbolic plane, and $X_{v_p}$ is a
$(p+1)$-valent regular tree.

If $\phi : \mathbf{PGL_2}(\mathbb{Z}[1/m]) \rightarrow
\mathbf{PGL_2}(\mathbb{Z}[1/m])$ is a quasi-isometry, then by
Lemma~\ref{l:1.4} we can replace $\phi$ by a quasi-isometry
$\phi_S$ which is the product of quasi-isometries $$\phi _\infty :
X_{v_\infty} \rightarrow X_{v_\infty}$$ and $$\phi _p :X_{v_p}
\rightarrow X_{v_{\tau (p)}}$$ for some permutation $\tau$ of the
primes dividing $m$.

If $m$ is prime, then this theorem reduces to Taback's theorem
\cite{Ta}. Now suppose $\ell$ is a prime dividing $m$. We let
$S'=\{\,v_\infty \,, v_\ell\,\}$ and $S''=\{\,v_\infty \,,
v_{\tau(\ell )}\,\}$.

It follows from the density of $\mathbf{PGL_2}(\mathbb{Z}[1/m])$
in
$$\prod _{p|m \, ; \, p \neq \tau ( \ell )} \mathbf{PGL_2}(\mathbb{Q}_p)$$
 that any point in $X_{S-S''}$ is a uniformly bounded distance from the
orbit $\mathbf{PGL_2}(\mathbb{Z}[1/m])x_{S-S''}$. Therefore, we
may assume that there is some $y \in
\mathbf{PGL_2}(\mathbb{Z}[1/m])x_{S-S''}$ such that $\phi_S \big(
X_{S'} \times \{x_{S-S'}\} \big) \subseteq X_{S''} \times \{y\}$.

Since $\phi(\Omega _S ) \subseteq \Omega _S$, we may assume that
$\phi_S \big( \Omega_S \cap ( X_{S'} \times \{x_{S-S'}\} ) \big)
\subseteq \Omega_S \cap ( X_{S''} \times \{y\} )$. It follows from
Lemma~\ref{l:1.3} that $\phi _S$ restricts to a quasi-isometry
between $\Omega _{S'}\times\{x_{S-S'}\}$ and $\Omega_{S''} \times
\{y\}$.

Note that by the product structure of $\phi _S$, we can assume
that $\phi _\infty \times \phi _ \ell$ restricts to a
quasi-isometry between $\Omega _{S'}$ and $\Omega _{S''}$ --- or
by Lemma~\ref{l:1.1}, a quasi-isometry between
$\mathbf{PGL_2}(\mathbb{Z}[1/\ell])$ and
$\mathbf{PGL_2}(\mathbb{Z}[1/{\tau (\ell)}])$. Taback showed that
for any such quasi-isometry we must have $\ell = \tau (\ell)$ and
that $\phi_\infty \times \phi _\ell$ is equivalent to a
commensurator $g_\infty \times g_\ell \in
\mathbf{PGL_2}(\mathbb{R}) \times \mathbf{PGL_2}(\mathbb{Q}_\ell)$
where $g_\infty$ and $g_\ell$ must necessarily be included in, and
represent the same element of, $\mathbf{PGL_2}(\mathbb{Q})$
\cite{Ta}.

As the above paragraph is independent of the prime $\ell|m$, the
element $g_\infty \in \mathbf{PGL_2}(\mathbb{Q})$ determines $\phi
_S$ and thus $\phi$.

\end{proof}

\bigskip

Having concluded our proof of the above special case, we return to
the general proof. Our goal is to show that $\phi _S$ is
equivalent to an element of
$\text{Comm}_{\text{Aut}(G_S)}(\mathbf{G}(\mathcal{O}_S))$. At
this point, the proof breaks into two cases.

\end{subsection}

\begin{subsection}{Case 1: $G_{V_K^\infty}$ is not locally
isomorphic to $\mathbf{PGL_2}(\mathbb{R})$}\label{s:2.5}

Notice that $\text{Comm}_{\text{Aut}(G_S)}(\mathbf{G}(\mathcal{O}_S))$ acts by
isometries on $X_S$. So a good first step toward our goal is to
show that $\phi _S$ is equivalent to an isometry. First, we will
show that $\phi _{V_K^\infty}$ is equivalent to an isometry.

\begin{lemma}\label{l:1.5} The quasi-isometry
$\phi _{V_K^\infty} : X_{V_K^\infty} \rt X_{V_K^\infty}$ is
equivalent to an isometry of the symmetric space $X_{V_K^\infty}$.
Indeed, it is equivalent to an element of
$\emph{Comm}_{\emph{Aut}(G_S)}(\mathbf{G}(\mathcal{O}))$.

\end{lemma}

\begin{proof} Notice that $\phi _{V_K^\infty}$ is
simply the restriction of $\phi _S$ to $X_{V_K^\infty} \times \{x
 _{S-V_K^\infty}\}$.

 Since $\mathbf{G}(\mathcal{O}_S)$ is
dense in $G_{S-V_K^\infty}$, the Hausdorff distance between
$\mathbf{G}(\mathcal{O_S})x_{S-V_K^\infty}$ and $X_{S-V_K^\infty}$
is finite. Thus, by replacing $\phi _S$ with an equivalent
quasi-isometry, we may assume that $X_{V_K^\infty} \times
\{x_{S-V_K^\infty}\}$ is mapped by $\phi _S$ into a space $X
_{V_K^\infty} \times \{y\}$ for some $y \in X_{S-V_K^\infty}$ with
$y \in \mathbf{G}(\mathcal{O}_S)x_{S-V_K^\infty}$.

  Since the Hausdorff distance between $\phi_S (\Omega _S)$ and $\Omega _S$ is finite,
   we have by Lemmas~\ref{l:1.3} and~\ref{l:1.1} that $\phi
_{V_K^\infty}$ induces a quasi-isometry of
$\mathbf{G}(\mathcal{O})$. The lemma follows from the existing
quasi-isometric classification of arithmetic lattices using our
assumption in this Case that $G_{V_K^\infty}$ is not locally
isomorphic to $\mathbf{PGL_2}(\mathbb{R})$; see \cite{Fa}.

\end{proof}

At this point, it is not difficult to see that Theorem~\ref{t:a} holds in
the case when every nonarchimedean factor of $G_S$ is higher rank:

\begin{lemma}\label{l:1.6} If $\emph{rank}_{K_v}(\mathbf{G})\geq2$
for all $v \in S-V_K^\infty$ then $\phi _S$ is equivalent to an
element of
$\emph{Comm}_{\emph{Aut}(G_S)}(\mathbf{G}(\mathcal{O}_S))$.
\end{lemma}

\begin{proof} By Lemma~\ref{l:1.4}, $\phi _v$ is equivalent to an isometry for all $v \in
S-V_K^\infty$. Combined with Lemma~\ref{l:1.5}, we know that $\phi _S$ is
equivalent to an isometry.

That $\phi_S$ is equivalent to an element of
$\text{Comm}_{\text{Aut}(G_S)}(\mathbf{G}(\mathcal{O}_S))$ follows from Proposition
7.2 of \cite{W2}. Indeed, any isometry of $X_S$ that preserves
$\Omega _S$ up to finite Hausdorff distance corresponds in a
natural way to an automorphism of $G_S$ that preserves
$\mathbf{G}(\mathcal{O}_S)$ up to finite Hausdorff distance, and
any such automorphism of $G_S$ is shown in Proposition 7.2 of
\cite{W2} to be a commensurator.

\end{proof}

For the remainder of Case 1, we are left to assume that there is
at least one $w \in S-V_K^\infty$ such that $\text{rank}_{K_w}
(\mathbf{G})=1$.

Before beginning the proof of the next and final lemma for Case 1,
it will be best to recall some standard facts about boundaries.

\begin{remark trees} If $w$ is a
nonarchimedean valuation of $K$, and
$\text{rank}_{K_w}(\mathbf{G})=1$, then $X_w$ is a tree.

For any minimal $K_w$-parabolic subgroup of $\mathbf{G}$, say
$\mathbf{M}$, we let $\varepsilon _\mathbf{M}$ be the end of $X_w$
such that $\mathbf{M}(K_w) \varepsilon _\mathbf{M}=\varepsilon
_\mathbf{M}$. Notice that the space of all ends of the form
$\varepsilon _\mathbf{P}$ where $\mathbf{P}$ is a minimal
$K$-parabolic subgroup of $\mathbf{G}$ forms a dense subset of the
space of ends of $X_w$.

\end{remark trees}

\begin{remark tits} For any
minimal $K$-parabolic subgroup of $\mathbf{G}$, say $\mathbf{P}$,
 we let $\delta _\mathbf{P}$ be the simplex in the
Tits boundary of $X_{V_K^\infty}$ corresponding to the group
$\prod _{v \in V_K^\infty}  { \mathbf{P}(K_v)}$.

If $\delta$ is a simplex in the Tits boundary of $X_{V_K^\infty}$,
and $\varepsilon$ is an end of the tree $X_w$, then we denote the
join of $\delta$ and $\varepsilon$ by $\delta * \varepsilon$. It
is a simplex in the Tits boundary of $X_T$ where $T=V_K^\infty
\cup \{w\}$.

\end{remark tits}

\begin{lemma}\label{l:1.7} Let $w \in S-V_K^\infty$ be such that
$\emph{rank}_{K_w} (\mathbf{G})=1$. Then $\phi _w : X_w
\rightarrow X_{\tau (w)}$ is equivalent to an isometry that is
induced by an isomorphism of topological groups $\mathbf{G}(K_w)
\rt \mathbf{G}(K_{\tau (w)})$.
\end{lemma}

\begin{proof} Below, we will denote the set of
valuations $V_K^\infty \cup \{\tau(w)\}$ by $T^\tau$.

We begin by choosing a minimal $K$-parabolic subgroup of
$\mathbf{G}$, say $\mathbf{P}$, and a geodesic ray $\rho :
\mathbb{R}_{\geq 0} \rt X_T$ that limits to the interior of the
simplex $\delta _\mathbf{P}
* \varepsilon _\mathbf{P}$.

By Lemma~\ref{l:1.5}, the image of $\phi _T \circ \rho$ under the
projection $X_T \rt X_{V_K^\infty}$ limits to a point in the
interior of $\delta _\mathbf{Q}$ for some minimal $K$-parabolic
subgroup $\mathbf{Q}$. Similarly, $\phi _w$ is a quasi-isometry of
a tree, so it maps each geodesic ray into a bounded neighborhood
of a geodesic ray that is unique up to finite Hausdorff distance.
Thus, the image of $\phi _T \circ \rho$ under the projection $X_T
\rt X_{\tau(w)}$ limits to $\varepsilon _\mathbf{Q'}$ for some
minimal $K_{\tau(w)}$-parabolic subgroup $\mathbf{Q'}$. Together,
these results imply that $\phi _T \circ \rho$ is a finite
Hausdorff distance from a geodesic ray that limits to a point in
the interior of $\delta _\mathbf{Q}*\varepsilon _\mathbf{Q'}$.

By Lemma~\ref{l:2.4}, there is a subspace $\mathcal{H}_\mathbf{P}$ of $X_T
$ corresponding to $\mathbf{P}$ (called a ``$T$-horoball'') such
that
$$t \mapsto d \Big(\rho (t) \,,\, X_T -
\mathcal{H}_\mathbf{P}\Big)$$ is an unbounded function. Thus,
$$t \mapsto d \Big(\phi _T \circ \rho (t) \,,\, X_{T^\tau} -
\phi _T (\mathcal{H}_\mathbf{P})\Big)$$ is also unbounded.

Using Lemma~\ref{l:2.1} and the fact that $\phi _T (\Omega _T)$ is a
finite Hausdorff distance from $\Omega _{T ^ \tau}$, we may
replace $\phi _T$ with an equivalent quasi-isometry to deduce that
$\phi _T(\mathcal{H}_\mathbf{P})$ is contained in the union of all
$T^\tau$-horoballs in $X_{T^\tau}$.

It will be clear from the definition given in Section~\ref{s:2}
that each $T$-horoball is connected. In addition,
Lemma~\ref{l:2.2} states that the collection of $T^\tau$-horoballs
is pairwise disjoint, so it follows that $\phi _T
(\mathcal{H}_\mathbf{P})$ is a finite distance from a single
$T^\tau$-horoball $\mathcal{H}_\mathbf{M} \se X_{T^\tau}$ where
$\mathbf{M}$ is a minimal $K$-parabolic subgroup of $\mathbf{G}$.
Therefore,
$$t \mapsto d \Big(\phi _T \circ \rho (t) \,,\, X_{T^\tau} -
\mathcal{H}_\mathbf{M}\Big)$$ is unbounded.

Because the above holds for all $\rho$ limiting to $\delta
_\mathbf{P}
* \varepsilon _\mathbf{P}$, and because $\phi _T \circ \rho$ limits to
 $\delta _\mathbf{Q}* \varepsilon _\mathbf{Q'}$, we have by Lemma~\ref{l:2.5},
that $\mathbf{Q}=\mathbf{M}=\mathbf{Q'}$. That is, $\phi
_{V_K^\infty}$ completely determines the map that $\phi _w$
induces between the ends of the trees $X_w$ and $X_{\tau(w)}$ that
correspond to $K$-parabolic subgroups of $\mathbf{G}$.

\bigskip \bigskip

Recall that by Lemma~\ref{l:1.5}, $\phi_{V_K^\infty}$ is
equivalent to a commensurator of $\mathbf{G}(\mathcal{O}) \leq
G_{V_K^\infty}$. Using Lemma 7.3 of \cite{W2}, $\phi_{V_K^\infty}$
(regarded as an automorphism of $G_{V_K^\infty}$) restricts to
$\mathbf{G}(K)$ as a composition
$$\delta \circ \sigma ^\circ : \mathbf{G}(K) \rt \mathbf{G}(K)$$
 where $\sigma$ is an automorphism of
$K$,
$$\sigma ^\circ :\mathbf{G}(K) \rt {^\sigma \mathbf{G}}(K)$$ is the
map obtained by applying $\sigma$ to the entries of elements in
$\mathbf{G}(K)$, and
$$\delta : {^\sigma \mathbf{G}} \rt \mathbf{G}$$ is a
$K$-isomorphism of $K$-groups.

Thus, if $\partial X_w$ and $\partial X_{\tau(w)}$ are the space
of ends of the trees $X_w$ and $X_{\tau(w)}$ respectively, and if
$\partial \phi _w :\partial X_w \rt \partial X_{\tau (w)}$ is the
boundary map induced by $\phi _w$, then we have shown that
$$\partial \phi _w (\varepsilon _\mathbf{P})=\varepsilon _{\delta
({^\sigma \mathbf{P}})}$$ for any $\mathbf{P} \leq \mathbf{G}$
that is a minimal $K$-parabolic subgroup of $\mathbf{G}$.

\bigskip \bigskip

Our next goal is to show that the valuation $w \circ \sigma ^{-1}$
is equivalent to $\tau(w)$. If this is the case, then $\delta
\circ \sigma ^\circ$ extends from a group automorphism of
$\mathbf{G}(K)$ to a topological group isomorphism
$$\alpha _w :\mathbf{G}(K_w)\rt \mathbf{G}(K_{\tau(w)})$$
If $\partial \alpha _w : \partial X_w \rt \partial X_{\tau(w)}$ is
the map induced by $\alpha _w$, then $\partial \alpha _w$ equals
$\partial \phi _w$ on the subset of ends in $\partial X_w$
corresponding to $K$-parabolic subgroups of $\mathbf{G}$ since
$\alpha _w$ extends $\delta \circ \sigma ^\circ$. Therefore,
$\partial \alpha _w =
\partial \phi _w$ on all of $\partial X_w$ by the density of the
 ``$K$-rational ends" in $\partial X_w$.
Thus, $\alpha _w$ determines $\phi _w$ up to equivalence. This
would prove our lemma.

So to finish the proof of this lemma, we will show that $w \circ
\sigma ^{-1}$ is equivalent to $\tau(w)$.

\bigskip \bigskip

For any maximal $K$-split torus $\mathbf{S}\leq \mathbf{G}$, we
let $\gamma ^w _\mathbf{S} \se X_w$ (resp. $\gamma
_\mathbf{S}^{\tau (w)} \se X_{\tau(w)}$) be the geodesic that
$\mathbf{S}(K_w)$ (resp. $\mathbf{S}(K_{\tau(w)})$) acts on by
translations.

Fix $\mathbf{S}$ and $\mathbf{T}$, two maximal $K$-split tori in
$\mathbf{G}$ such that $\gamma ^w _\mathbf{S} \cap \gamma ^w
_\mathbf{T}$ is nonempty and bounded. We choose a point $a \in
\gamma ^w _\mathbf{S} \cap \gamma ^w _\mathbf{T}$.

Since $\mathbf{S}(\mathcal{O}_T)$ is dense in $\mathbf{S}(K_w)$,
there exists a group element $g_n \in \mathbf{S}(\mathcal{O}_T)$
for each $n \in \mathbb{N}$ such that
$$d \Big( g_n (\gamma ^w_\mathbf{T})\,,\,a\Big)>n$$
Note that $g_n (\gamma ^w_\mathbf{T}) = \gamma ^w
_{g_n\mathbf{T}g_n^{-1}}$. Thus
$$d \Big( \phi _w (\gamma ^w
_{g_n\mathbf{T}g_n^{-1}}) \,,\, \phi _w(a)\Big) $$ is an unbounded
sequence.

As ${g_n\mathbf{T}g_n^{-1}}$ is $K$-split, $\phi _w(\gamma ^w
_{g_n\mathbf{T}g_n^{-1}})$ is a uniformly bounded Hausdorff
distance from
$$\gamma ^{\tau(w)} _{\delta \circ \sigma ^\circ
({g_n\mathbf{T}g_n^{-1}})}=\delta \circ \sigma ^\circ (g_n) \gamma
^{\tau(w)} _{\delta \circ \sigma ^\circ (\mathbf{T})}$$ because a
geodesic in $X_{\tau(w)}$ is determined by its two ends.

We finally have that
$$d \Big(\delta \circ \sigma ^\circ (g_n) \gamma ^{\tau(w)} _{\delta
\circ \sigma ^\circ (\mathbf{T})} \,,\, \phi _w (a)\Big)$$ is an
unbounded sequence. It is this statement that we shall contradict
by assuming that $w \circ \sigma ^{-1}$ is inequivalent to
$\tau(w)$.

Note that $g_n \in \mathbf{G}(\mathcal{O}_v)$ for all $v \in
V_K-T$ since $g_n \in \mathbf{S}(\mathcal{O}_T)$. Thus, $\sigma ^0
(g_n) \in {^\sigma \mathbf{G}}(\mathcal{O}_{v\circ \sigma ^{-1}})$
for all $v \in V_K-T$. If it were the case that $w \circ \sigma
^{-1}$ is inequivalent to $\tau(w)$, then it follows that $\sigma
^0 (g_n) \in {^\sigma \mathbf{G}}(\mathcal{O}_{\tau(w)})$. Hence,
$\delta \circ \sigma ^0 (g_n)$ defines a bounded sequence in
$\mathbf{G}(K_{\tau(w)})$. Therefore, $$d \Big(\delta \circ \sigma
^\circ (g_n) \gamma ^{\tau(w)} _{\delta \circ \sigma ^\circ
(\mathbf{T})}\,,\, \phi _w (a)\Big)$$ is a bounded sequence, our
contradiction.

\end{proof}

The proof of Theorem~\ref{t:a} in Case 1 is complete with the
observation that applications of Lemma~\ref{l:1.7} to tree
factors, allows us to apply the Proposition 7.2 of \cite{W2} as we
did in Lemma~\ref{l:1.6}.

\end{subsection}

\begin{subsection}{Case 2: $G_{V_K^\infty}$ is locally
isomorphic to $\mathbf{PGL_2}(\mathbb{R})$}

It follows that $V_K^\infty$ contains a single valuation $v$, and
that $K_v \cong \mathbb{R}$. Thus $K=\mathbb{Q}$, and $V_K^\infty$
is the set containing only the standard real metric on
$\mathbb{Q}$.

Our assumption that $\mathbf{G}$ is absolutely simple implies that
$G_{V_K^\infty}$ is actually isomorphic to
$\mathbf{PGL_2}(\mathbb{R})$. Thus, $\mathbf{G}$ is a
$\mathbb{Q}$-form of $\mathbf{PGL_2}$. As we are assuming that
$\mathbf{G}$ is $\mathbb{Q}$-isotropic, it follows from the
classification of $\mathbb{Q}$-forms of $\mathbf{PGL_2}$ that
$\mathbf{G}$ and $\mathbf{PGL_2}$ are $\mathbb{Q}$-isomorphic (see
e.g. page 55 of \cite{Ti class}).

From our assumptions in the statement of Theorem~\ref{t:a}, $S
\neq V_K^\infty$. As the only valuations, up to scale, on
$\mathbb{Q}$ are the real valuation and the $p$-adic valuations,
$\mathbf{G}(\mathcal{O}_S)$ is commensurable with
$\mathbf{PGL_2}(\mathbb{Z}[1/m])$ for some $m \in \mathbb{N}$ with
$m \neq 1$. Thus, Case 2 of Theorem~\ref{t:a} follows from
Theorem~\ref{t:m}.

Our proof of Theorem~\ref{t:a} is complete modulo the material
from Section~\ref{s:2}.

\end{subsection}

\end{section}

\begin{section}{Horoball patterns
 in a product of \newline \indent $\,$ a tree and a symmetric space}\label{s:2}

In this section we will study the components of $X_S-\Omega_S$
when $X_S$ is a product of a symmetric space and a tree.

\begin{remark note} We let $w$ be a
 nonarchimedean valuation on
$K$ such that $\text{rank}_{K_w}(\mathbf{G})=1$. Then we set $T$
equal to $V_K^\infty \cup \{w \}$.

\end{remark note}

\begin{subsection}{Horoballs in rank one symmetric spaces.}
Let $\mathbf{P}$ be a minimal $K$-parabolic subgroup of
$\mathbf{G}$. As in the previous section, we let $\delta
_\mathbf{P}$ be the simplex in the Tits boundary of
$X_{V_K^\infty}$ corresponding to the group $\prod _{v \in
V_K^\infty}  \mathbf{P}(K_v)$.

Note that $\mathbf{G}$ being $K$-isotropic and
$\text{rank}_{K_w}(\mathbf{G})=1$ together implies that
$\text{rank}_{K}(\mathbf{G})=1$. Borel proved that the latter
equality implies that $X_{V_K^\infty} - \Omega _{V_K^\infty}$ can
be taken to be a disjoint collection of horoballs (17.10
\cite{Bo}).

To any horoball of $X_{V_K^\infty} - \Omega _{V_K^\infty}$, say
$H$, there corresponds a unique $\delta _ \mathbf{P}$ as above
such that any geodesic ray $\rho : \mathbb{R}_{\geq 0} \rt
X_{V_K^\infty}$ that limits to $\delta _\mathbf{P}$ defines an
unbounded function
$$t \mapsto d\Big( \rho (t) \,,\, X_{V_K^\infty} - H \Big)$$

\end{subsection}

\begin{subsection}{$T$-horoballs in $X_T$} Let $y \in X_w$ and suppose
$y \in \mathbf{G}(\mathcal{O}_T) x_w$. Recall that by Lemma~\ref{l:1.3},
the space $\Omega _T \cap \big(X_{V_K^\infty} \times \{y\} \big)$
is a finite Hausdorff distance from $\Omega_{V_K^\infty} \times
\{y\}$.

For any minimal $K$-parabolic subgroup of $\mathbf{G}$, say
$\mathbf{P}$, we let $\mathcal{H}_{\mathbf{P},\infty}^y \subseteq
X_{V_K^\infty} \times \{y\}$ be the horoball of $\Omega _T \cap
\big(X_{V_K^\infty} \times \{y\} \big)$ that corresponds to
$\delta _\mathbf{P}$.

For arbitrary $x \in X_w$, we define
$$\mathcal{H}_{\mathbf{P},\infty}^x =
\mathcal{H}_{\mathbf{P},\infty}^y$$ where $y \in
\mathbf{G}(\mathcal{O}_T) x_w $ minimizes the distance between $x$
and $\mathbf{G}(\mathcal{O}_T) x_w$.

We let
$$\mathcal{H}_\mathbf{P} = \bigcup _{x \in X_w} \Big(
\mathcal{H}_{\mathbf{P},\infty}^x \times \{ x \} \Big)$$ Each of
the spaces $\mathcal{H}_\mathbf{P}$ is called a
\emph{$T$-horoball}.

Let $\mathcal{P}$ be the set of all minimal $K$-parabolic
subgroups of $\mathbf{G}$. The following lemma follows directly
from our definitions. It will be used in the proof of
Lemma~\ref{l:1.7}.

\begin{lemma}\label{l:2.1} The Hausdorff distance between $X_T-\Omega_T$
and
$$\bigcup _{\mathbf{P} \in \mathcal{P}} \mathcal{H}_\mathbf{P}$$ is finite.
\end{lemma}

We record another observation to be used in the proof of Lemma~\ref{l:1.7}.

\begin{lemma}\label{l:2.2}  If $\mathbf{P} \neq \mathbf{Q}$ are minimal
$K$-parabolic subgroups of $\mathbf{G}$, then
$\mathcal{H}_\mathbf{P} \cap \mathcal{H}_\mathbf{Q} = \emptyset$.
\end{lemma}

\begin{proof} The horoballs comprising $\big(
X_{V_K^\infty} - \Omega _{V_K^\infty} \big) \times \{x_w\}$ are
pairwise disjoint, and are a finite Hausdorff distance from the
horoballs of $\Omega _T \cap \big(X_{V_K^\infty} \times \{x_w\}
\big)$ by Lemma~\ref{l:1.2}. Hence, if $y = g x_w$ for some $g \in
\mathbf{G}(\mathcal{O}_T)$, then the horoballs determined by $$
\Omega _T \cap \big(X_{V_K^\infty} \times \{y\} \big)  = g \Big[
\Omega _T \cap \big(X_{V_K^\infty} \times \{x_w\} \big) \Big]
$$ are disjoint.

\end{proof}

\end{subsection}

\begin{subsection}{Deformations of horoballs over geodesics in $X_w$}
We let $\pi : X_T \rt X_{V_K^\infty}$ be the projection map. Note
that if $x \in X_w$ and $\mathbf{P}$ is a minimal $K$-parabolic
subgroup of $\mathbf{G}$, then $\pi
(\mathcal{H}^x_{\mathbf{P},\infty})$ is a horoball in
$X_{V_K^\infty}$ that is based at $\delta _\mathbf{P}$.

Recall that for any minimal $K_w$-parabolic subgroup of
$\mathbf{G}$, say $\mathbf{Q}$, we denote the point in the
boundary of the tree $X_w$ that corresponds to $\mathbf{Q}(K_w)$
by $\varepsilon _\mathbf{Q}$.

\begin{lemma}\label{l:2.3} Suppose $ \mathbf{P}$ is a minimal
$K$-parabolic subgroup of $\mathbf{G}$ and that $ \mathbf{Q}$ is a
minimal $K_w$-parabolic subgroup of $\mathbf{G}$. If $\gamma :
\mathbb{R} \rightarrow X_w$ is a geodesic with $\gamma (\infty )
=\varepsilon _\mathbf{P}$ and $\gamma ( - \infty ) =\varepsilon
_\mathbf{Q }$ then
\begin{quote} (i) $s \leq t$ implies $\pi \big(\mathcal{H}_{\mathbf{P},\infty}^{\gamma
(s)}\big)
 \subseteq \pi \big( \mathcal{H}_{\mathbf{P},\infty}^{\gamma
(t)}\big)$ \\
\medskip (ii) $\cup _{t\in \mathbb{R}} \pi \big(
\mathcal{H}_{\mathbf{P},\infty}^{\gamma
(t)} \big) = X_{V_K^\infty}$ \\
\medskip (iii) $\cap _{t\in \mathbb{R}}
\big( \mathcal{H}_{\mathbf{P},\infty}^{\gamma (t)} \big) = \emptyset$ \\
\medskip (iv) There exists constants $L_\mathbf{P},C>0$ such that if $$h(s,t)=
d \Big( \pi \big(\mathcal{H}_{\mathbf{P},\infty}^{\gamma (s)}
\big)\,,\,\pi \big(\mathcal{H}_{\mathbf{P},\infty}^{\gamma (t)}
\big)\Big)$$ \hspace*{.2in}then
$$|\,h(s,t)-L_\mathbf{P}|s-t|\,|\leq C$$
\end{quote}

\end{lemma}

\begin{proof} As the ends of $X_w$ corresponding to
$K$-parabolic subgroups are a dense subset of the full space of
ends, it suffices to prove this lemma when $\mathbf{Q}$ is defined
over $K$. In this case, the image of $\gamma $ corresponds to a
$K$-split torus $\mathbf{S} \leq \mathbf{G}$ that is contained in
$\mathbf{P}$.

Let $\alpha$ be a root of $\mathbf{G}$ with respect to
$\mathbf{S}$ such that $\alpha$ is positive in $\mathbf{P}$. Since
the diagonal embedding of $\mathbf{S}(\mathcal{O}_T)$ in the group
$\prod_{v\in V_K^\infty}\mathbf{S}(K_v)$ has a dense image, there
is some $b \in \mathbf{S}(\mathcal{O}_T)$ such that $|\alpha
(\mathbf{S}(b))|_v <1$ for all $v \in V_K^\infty$. Thus,
 $\mathbf{S}(b) \pi \big(\mathcal{H}_{\mathbf{P},\infty}^{\gamma (0)}\big)$ is a horoball
in $X_{V_K^\infty}$ that strictly contains $\pi \big(
\mathcal{H}_{\mathbf{P},\infty}^{\gamma (0)}\big)$. Generally, we
have
$$\mathbf{S}(b)^m \pi \big( \mathcal{H}_{\mathbf{P},\infty}^{\gamma (0)}\big) \subsetneq
 \mathbf{S}(b)^n \pi \big(
\mathcal{H}_{\mathbf{P},\infty}^{\gamma (0)}\big)$$ for all $m,n
\in \mathbb{Z}$ with $m<n$.

By the product formula we have $|\alpha (\mathbf{S}(b))|_w >1$.
Thus, there is a positive number $\lambda >0$ such that
 $\gamma (n \lambda ) = \mathbf{S}(b)^n \gamma (0)$ for any
  $n \in \mathbb{Z}$. It follows for $m<n$ that
$$\pi \big( \mathcal{H}_{\mathbf{P},\infty}^{\gamma
(m\lambda)}\big) = \mathbf{S}(b)^m \pi \big(
\mathcal{H}_{\mathbf{P},\infty}^{\gamma (0)}\big) \subsetneq
\mathbf{S}(b)^n \pi \big( \mathcal{H}_{\mathbf{P},\infty}^{\gamma
(0)}\big) = \pi \big( \mathcal{H}_{\mathbf{P},\infty}^{\gamma
(n\lambda )}\big)$$

We let $$L_\mathbf{P} = \frac{1}{\lambda} d \Big(\pi
\big(\mathcal{H}_{\mathbf{P},\infty}^{\gamma (0)} \big) \,,\, \pi
\big(\mathcal{H}_{\mathbf{P},\infty}^{\gamma (\lambda)} \big)
\Big)$$ so that $$h(m \lambda ,n \lambda )=L_\mathbf{P} \lambda
|m-n|=L_\mathbf{P} |m \lambda - n \lambda |$$ Then we take
$$C' =\max_{0\leq s \leq t \leq \lambda}\Big\{ d \Big(\pi
\big(\mathcal{H}_{\mathbf{P},\infty}^{\gamma (s)} \big) \,,\, \pi
\big(\mathcal{H}_{\mathbf{P},\infty}^{\gamma (t)} \big) \Big)
\Big\}$$ and, say, $$C=2C'+L_\mathbf{P}d\big(\gamma
(0),\gamma(\lambda)\big)$$

\end{proof}

\end{subsection}

\begin{subsection}{Basepoints in the Tits boundary for $T$-horoballs} The
Tits boundary for $X$ is the spherical join of the Tits boundary
for the symmetric space $X_{V^\infty _K}$ and the Tits boundary
for the tree $X_w$.

The purpose of the following two lemmas---and of this entire
section---is to show that each $T$-horoball
$\mathcal{H}_\mathbf{P}$ is geometrically associated with the join
of $\delta _\mathbf{P}$ and $\varepsilon _\mathbf{P}$, denoted
$\delta _\mathbf{P} * \varepsilon _\mathbf{P}$.

\begin{lemma}\label{l:2.4} Let $ \mathbf{P}$ be a minimal
$K$-parabolic subgroup of $\mathbf{G}$. Any geodesic ray $\rho
:\mathbb{R}_{\geq 0} \rightarrow X_T$ that limits to the simplex
$\delta _\mathbf{P} * \varepsilon _\mathbf{P}$ in the Tits
boundary of $X_T$ defines an unbounded function when composed with
the distance from the complement of $\mathcal{H}_\mathbf{P} $ in
$X_T$:
$$t \mapsto d\Big( \rho (t) \,,\, X_T - \mathcal{H}_\mathbf{P} \Big)$$
\end{lemma}

\begin{proof} Any such geodesic ray $\rho$ is a
product of a geodesic ray
 $b: \mathbb{R}_{\geq 0} \rightarrow X_{V_K^\infty}$ that limits to
  $\delta _\mathbf{P}$ and a geodesic ray
$c: \mathbb{R}_{\geq 0} \rightarrow X _w$ that limits to
$\varepsilon _\mathbf{P}$.

Let $Y = \pi \big(\mathcal{H}_{\mathbf{P},\infty}^{c (0)} \big)
\times c(\mathbb{R}_{\geq 0})$. Since $$t \mapsto d\Big( b (t)
\,,\, X_{V_K^\infty} - \pi \big(\mathcal{H}_{\mathbf{P},\infty}^{c
(0)} \big) \Big)$$ is unbounded, $$t \mapsto d\Big( \rho (t) \,,\,
X_T - Y \Big)$$ is unbounded. The lemma follows from
Lemma~\ref{l:2.3}\emph{(i)} which guarantees that $Y \subseteq
\mathcal{H}_{\mathbf{P}}$.

\end{proof}

\begin{lemma}\label{l:2.5} Suppose
$\mathbf{Q}$ and $\mathbf{M}$ are
 minimal
$K$-parabolic subgroups of $\mathbf{G}$, and that $\mathbf{Q'}$ is
a minimal $K_w$-parabolic subgroup of $\mathbf{G}$ with
$\mathbf{M}\neq \mathbf{Q}$ or $\mathbf{M}\neq \mathbf{Q'}$. Then
there is a geodesic ray $\rho :\mathbb{R}_{\geq 0} \rightarrow X$
with $\rho (\infty) \in \delta _\mathbf{Q} * \varepsilon
_\mathbf{Q'}$ such that the function $$t \mapsto d\Big( \rho (t)
\,,\, X_T - \mathcal{H}_\mathbf{M} \Big)$$ is bounded.

\end{lemma}

\begin{proof} Choose a geodesic ray
 $b: \mathbb{R}_{\geq 0} \rightarrow X_{V_K^\infty}$ that limits
 to $\delta _\mathbf{Q}$ and a geodesic ray $c: \mathbb{R}_{\geq 0}
\rightarrow X _w$ that limits to $ \varepsilon _\mathbf{Q'}$. Let
$r$ be the ratio of the speed of $b$ to the speed of $c$.

If $\mathbf{M}\neq \mathbf{Q'}$, then after ignoring at most a
bounded interval of $c$, we can extend $c$ to a bi-infinite
geodesic with $c(-\infty) = \varepsilon _\mathbf{M}$. With
$L_\mathbf{M}$ as in Lemma~\ref{l:2.3}, $\rho(t) = \big( b(L_\mathbf{M} t)
\,,\, c( r  t) \big)$ defines a geodesic ray satisfying the lemma.

In the remaining case, $\mathbf{M}\neq \mathbf{Q}$ and $\mathbf{M}
= \mathbf{Q'}$.

The distance from $b(t)$ to $\pi
\big(\mathcal{H}^{b(0)}_{\mathbf{M},\infty}\big)$ is a convex
function in $t$. Since $\mathbf{M}\neq \mathbf{Q}$, this function
has a positive derivative, $u>0$, for some large value of $t$.
Then $\rho(t) =\big( b( L_\mathbf{M}t) \,,\, c(ur t)\big)$ defines
a geodesic ray satisfying the lemma.

\end{proof}

\end{subsection}

\end{section}

\bigskip

\noindent Kevin Wortman
\newline \noindent Department of Mathematics
\newline \noindent University of Utah
\newline \noindent 155 South 1400 East, Room 233
\newline \noindent Salt Lake City, UT 84112-0090
\newline \noindent wortman@math.utah.edu

\end{document}